\documentclass[11pt]{amsart}
\usepackage{amsmath, amsthm, amsfonts, amssymb, graphicx}
\usepackage{setspace}
\usepackage{fullpage}
\usepackage{amscd}
\usepackage{amsaddr}
\theoremstyle{plain} \numberwithin{equation}{section}

\newtheorem*{theorem*}{Theorem}

\def\00{(0,0)}
\def\11{(1,1)}

\begin{document}

\title{Gradient descent in some simple settings}
\author{Y.\ Cooper}
\email{yaim@math.ias.edu}
\date{\today}
\maketitle

\vspace{-.1in}

\begin{abstract}
In this note, we observe the behavior of gradient flow and discrete and noisy gradient descent in some simple settings.  It is commonly noted that addition of noise to gradient descent can affect the trajectory of gradient descent.  Here, we run some computer experiments for gradient descent on some simple functions, and observe this principle in some concrete examples.\end{abstract}

\section{Introduction}
In this note, we are interested in the behavior of discrete and noisy gradient descent in settings where there are minima with different characteristics.  To study this in a very simple setting, we consider gradient descent on two periodic functions.  The first we consider is 

$$f(x) = \sin(\pi x) + \cos(2 \pi x)+2.$$

This function is periodic, and has wells of two different depths, and slightly different widths as well.  The shallower wells have depth approximately 1 and width approximately 0.8, while the deeper wells have depth approximately 3 and width approximately 1.  

\begin{figure}[h]
\includegraphics[width=2.8in]{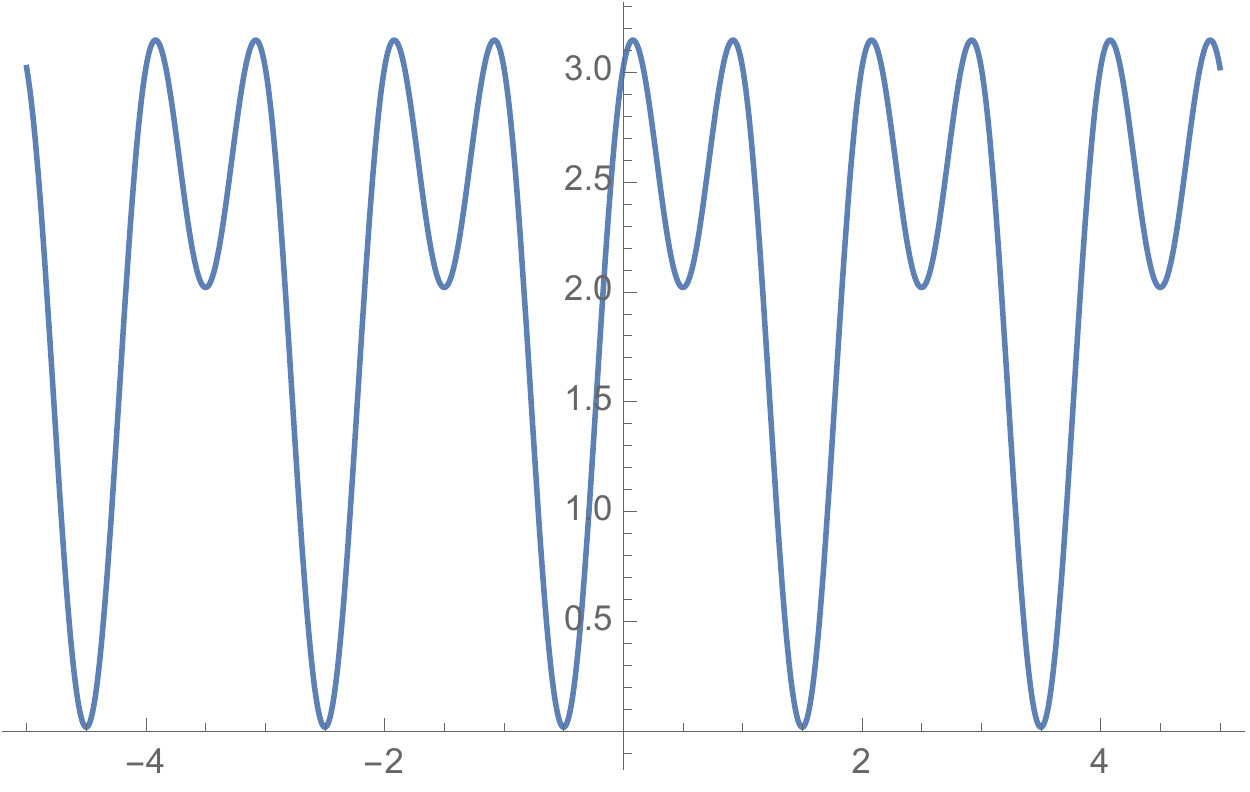}
\caption{We consider gradient descent on the simple function $f(x) = \sin(\pi x) + \cos(2 \pi x)+2$, which has wells of two different depths.}
\end{figure}

The second we consider is 
$$g(x) = \left(\sin(\pi x) + \frac{\sin(2 \pi x)}{2}\right)^2.$$

This function is periodic, and has two different types of wells.  This time they have the same depth, but are of different widths.  The narrower wells have width 2/3, while the wider ones have width 4/3.

\begin{figure}[h]
\includegraphics[width=3.5in]{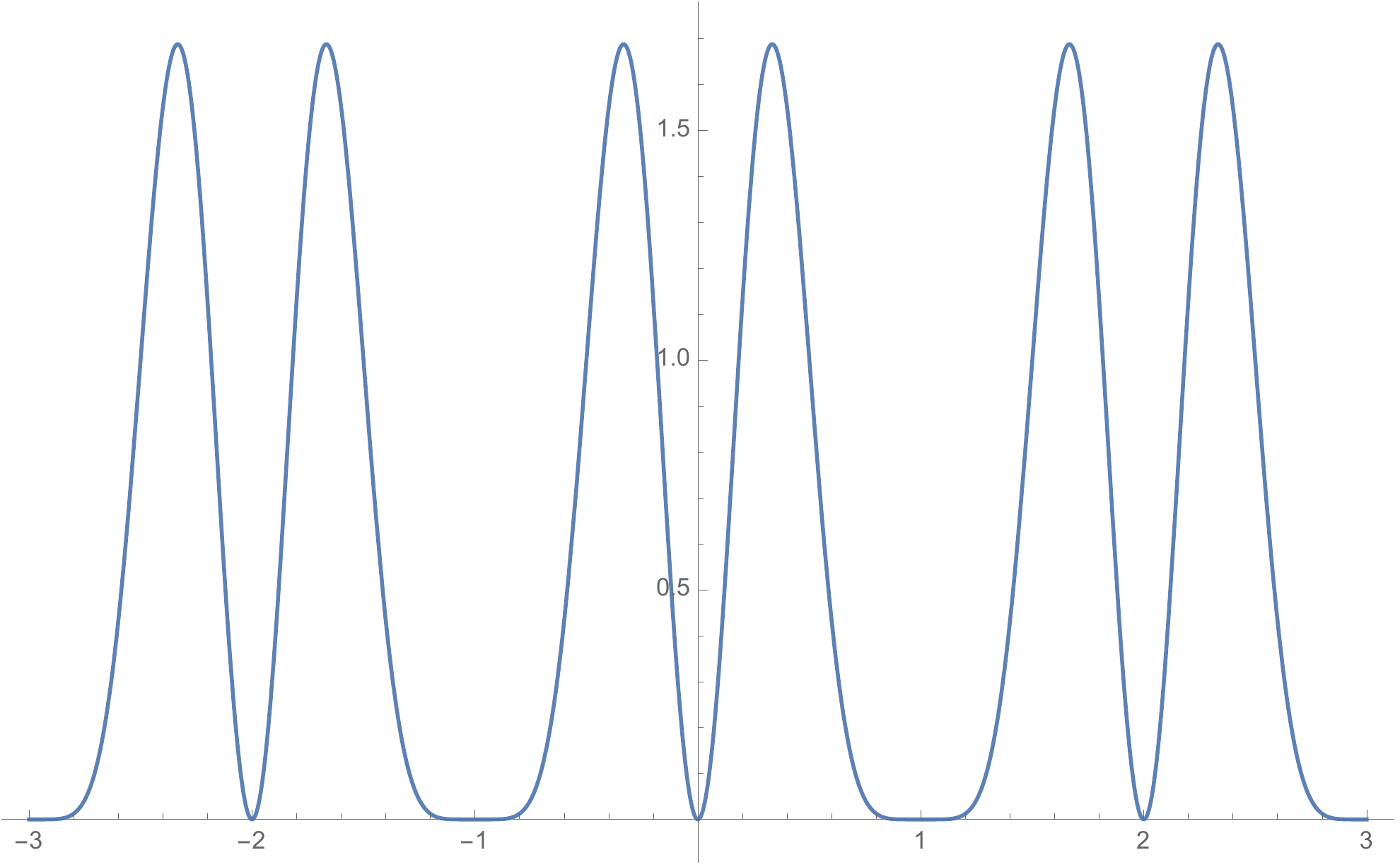}
\caption{We consider gradient descent on the simple function $g(x) = \left(\sin(\pi x) + \frac{\sin(2 \pi x)}{2}\right)^2$, which has wells of two different widths.}
\end{figure}

In this note, we explore the behavior of gradient flow and of noisy discrete gradient descent in these landscapes.

\subsection{Acknowledgements}

The author would like to thank Avi Wigderson for conversations that inspired this experiment and warm encouragement, and Nathaniel Bottman for help with the experiments.

\section{Gradient descent on $f(x)$: Theory}

In this section, we will consider gradient flow as well as a noisy form of discrete gradient descent.  In all cases, we will consider the setting where the initial position is drawn from a uniform distribution on a finite interval $[a,b]$.  

\subsection{Gradient flow}

First, we consider the behavior of gradient flow on the function $f(x) = \sin(\pi x) + \cos(2 \pi x)+2$ on an interval $[a,b]$, with a given starting point $p_0 \in [a,b]$.  For simplicity, we assume that $a,b$ are maxima of $f$, and that the interval $[a,b]$ contains the same number of shallow wells as deep wells.  

With measure zero, $p_0$ will be a critical point of $f$.  Assuming that $p_0$ is not a critical point of $f$, there are two possibilities.  Either $p_0$ is in the basin of attraction of a shallow well $S$, or it is in the basin of attraction of a deep well $D$.  While neither $S$ or $D$ is convex, nonetheless, under gradient flow if $p_0$ is in the basin of attraction of $S$ the flow line will end in the unique minimum of $S$, and similarly for $D$.  

Therefore, if we initialize $p_0$ randomly, we should expect that the ratio of the probability that a local nonglobal minimum is found to the probability that a global minimum is reached is the ratio of the width of the basin of attraction of $S$ to the width of the basin of attraction of $D$.  The width of the basin of attraction of $S$ is approximately 0.84, and the width of the basin of attraction of $D$ is approximately 1.16, so their ratio is approximately 0.72.  So we expect that under uniform random initialization, a minimum of a shallow well will be found approximately 0.72 as often as a minimum of a deep well.

\subsection{Discrete gradient descent with $\epsilon$-jitter}\label{with jitter} \label{jitterthy}

Having considered the case of gradient flow, which discrete gradient descent approximates, we now turn to modifications of discrete gradient descent which add some noise to the process.  Discrete gradient descent and noisy gradient descent have been used and studied by many, one reference for noise is \cite{AR1}.  One motivation for doing so is that perhaps in doing so one can bias the procedure toward deeper minima.  

One way to introduce noise into the process of discrete gradient descent is to add at each step a small random vector to the gradient vector.  We call this modification discrete gradient descent with $\epsilon$-jitter.  

In this case, we begin at some initial position $p_0$.  At the $t+1^{st}$ step, we let
$$
p_{t+1} = p_t - \tau \nabla L (p_t) - (\epsilon_{t})
$$
where $\epsilon_{t}$ are drawn from a gaussian distribution with norm 0 and standard deviation $\epsilon$.  

One consideration is the step size $\tau$.  A reasonable regime to consider is $0 < \tau < w/g$, where $w$ is the average width of wells and $g$ is the average magnitude of the gradient of $f$ in the interval under consideration.  This is because for $\tau$ larger than $w/g$, on most iterations the position $p_t$ will go from one well to another, and the process will become more like jumping between wells than like gradient descent.  

For step sizes in the regime $0 < \tau < w/g$, we expect that the addition of noise of order $\epsilon$ to the gradient descent algorithm will affect the relative probabilities of finding deep and shallow minima.  In the next section, we will experimentally study the dependence of that phenomenon on $\epsilon$.

\section{Discrete gradient descent on $f(x)$: Computer experiments}

Having established the processes we wish to study as well as our theoretical expectations of the parameter regimes in which they should be implemented and the expected behavior, we now turn to some computer experiments.  In all the experiments of this section, we run variants of gradient descent on the function $f(x) = \sin(\pi x) + \cos(2 \pi x)+2$ over the interval $[ -5.92, 6.08]$, and initialize from the uniform distribution on that interval.

\subsection{Gradient flow}

First, we approximate gradient flow by discrete gradient descent with small step size.  We find, as expected, that the ratio of the probabilities of arriving at the minimum of a shallow well versus a deep well is equal to the ratio of the basins of attraction.  

\begin{figure}[h]
\includegraphics[width=3in]{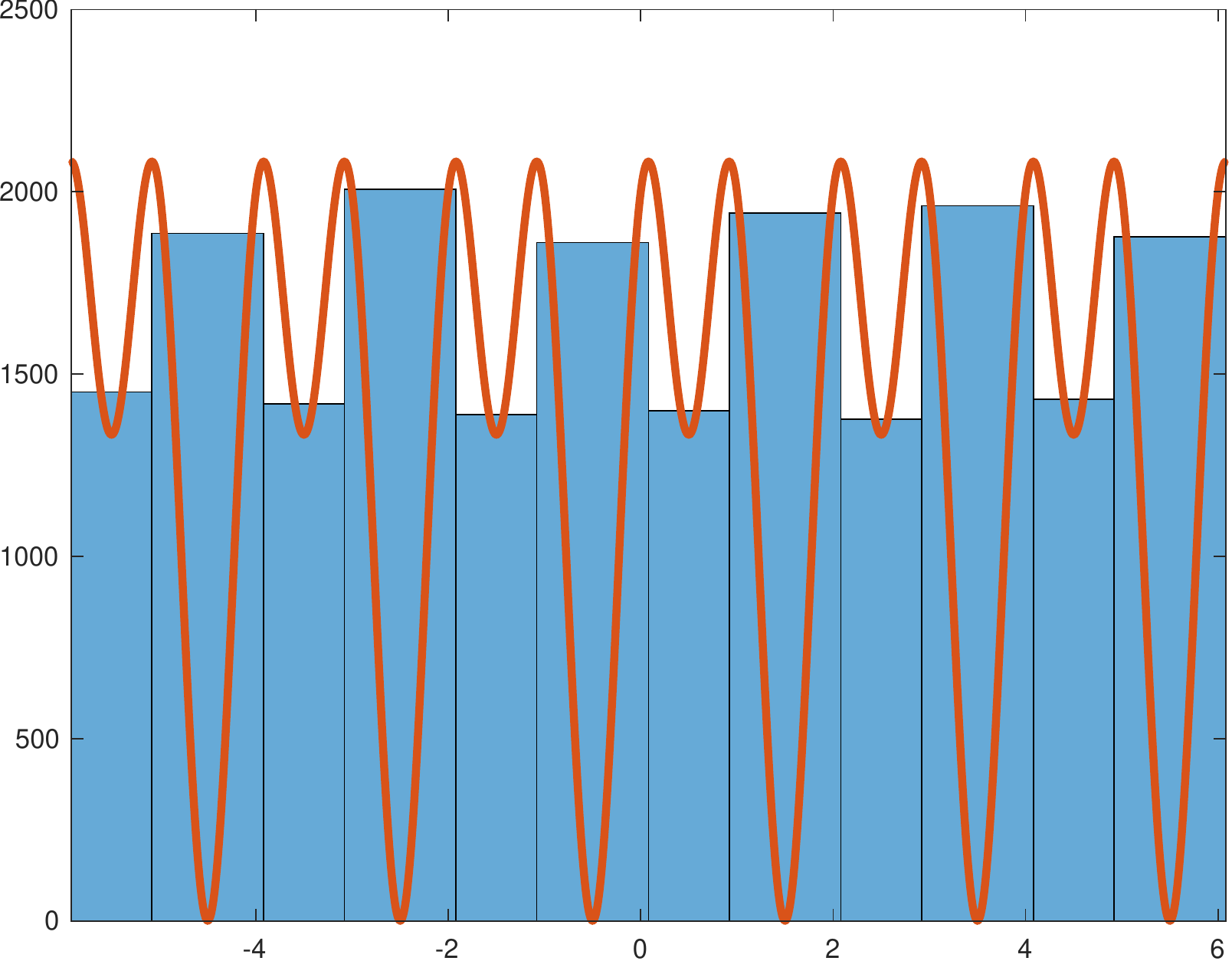}
\caption{In this experiment, we implement gradient descent with step size $\tau = 0.01$ on the function $f(x) = \sin(\pi x) + \cos(2 \pi x)+2$ over the interval $[ -5.92, 6.08]$.  We randomly initialize from the uniform distribution on this interval, and run the experiment 20,000 times.  Above is a histogram which shows the number of times the process ended in each well.}
\end{figure}

Namely, the ratio we find in the experiment described above is 0.73 , while the ratio between the widths of the basins of attraction of the two kinds of wells was computed above as 0.72.  

\subsection{Discrete gradient descent with $\epsilon$-jitter}\label{with jitter}

In this section, we outline the results of some computer experiments carrying out discrete gradient descent with $\epsilon$-jitter.  We find that indeed, for some values of $\epsilon$, adding noise to the gradient descent algorithm does bias the procedure toward finding the deeper minima over the shallower ones.

In the first set of experiments, we implement gradient descent with step size $\tau = .01$ on the function $f(x) = \sin(\pi x) + \cos(2 \pi x)+2$ over the interval $[ -5.92, 6.08]$.  We randomly initialize from the uniform distribution on this interval, and run the experiment 20,000 times.  We do this for several different values of noise $\epsilon$, and record the the number of times the process ended in each well.   

When $\epsilon=0$, the experiment went very similarly to the experiment approximating gradient flow, as expected.  With a relatively small amount of added noise however, the picture becomes dramatically different.

\begin{figure}[h]
\includegraphics[width=4in]{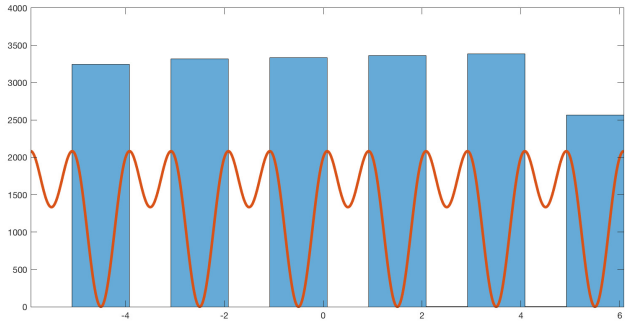}
\caption{This histogram shows the number of times the process ended in each well with $\epsilon = 0.15$.  The ratio of the probability of landing in a shallow well to the probability of landing in a deep well was 0.0003.}
\end{figure}

With a small amount of noise added to the process at each step, the behavior of gradient descent on this function $f$ becomes dramatically different.  With $\epsilon = 0.15$, the probability of finding a shallow minimum becomes nearly zero.  

\begin{figure}[h]
\includegraphics[width=4in]{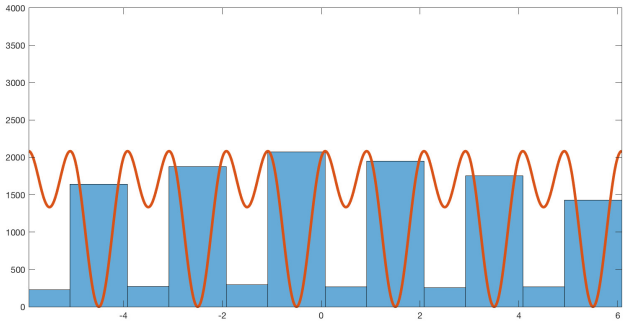}
\caption{This histogram shows the number of times the process ended in each well with $\epsilon = 0.25$.  The ratio of the probability of landing in a shallow well to the probability of landing in a deep well was 0.15.}
\end{figure}

The above behavior holds for a range of choices for $\epsilon$.  At $\epsilon = 0.25$, the effect is less dramatic, but the modified form of gradient descent continues to be strongly biased toward deep minima.

\begin{figure}[h]
\includegraphics[width=4in]{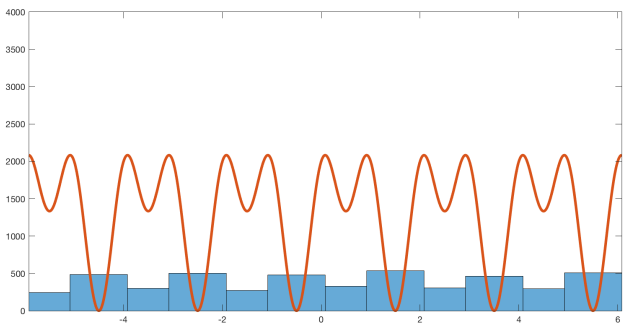}
\caption{This histogram shows the number of times the process ended in each well with $\epsilon = 0.5$.  The ratio of the probability of landing in a shallow well to the probability of landing in a deep well was 0.58.}
\end{figure}

Even for large values of $\epsilon$, when the noise is so large as to cause many trials to leave the interval of interest, the bias toward deep minima persists.  All the trials are initialized in the interval $[ -5.92, 6.08]$, but over the course of the 1,000 steps, most trials leave the interval.  Nonetheless, for those that remain in the interval, there is still a bias toward deep minima.

To address more systematically the question of the dependence of the behavior of discrete gradient descent with $\epsilon$-jitter on the step size $\tau$ and noise parameter $\epsilon$, we run experiments for a range of $\tau$ and $\epsilon$.  The results are collected in the following set of charts.  

\begin{figure}[h]\label{6}
\includegraphics[width=4.5in]{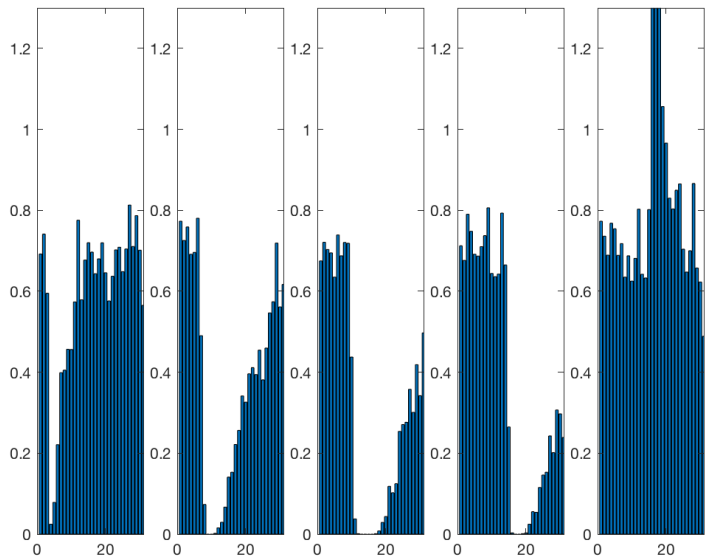}
\caption{These charts, from left to right, concern the behavior of noisy gradient descent on $f(x) = \sin(\pi x) + \cos(2 \pi x)+2$, with step size $\tau = 0.001, 0.01, 0.02, 0.04$, and $0.06$.  Each chart shows, for that step size, how the ratio of the probability of ending in a deep vs a shallow well varies as the magnitude $\epsilon$ of the added noise varies from $0$ to $0.5$.}
\end{figure}

Each chart shows the ratio $r$ of the probability of finding a shallow minimum to the probability of finding a deep minimum when the step size $\tau$ is fixed, and the noise $\epsilon$ goes from 0 to $0.5$.  That is, $\epsilon$ is on the $x$-axis, and the ratio $r$ is plotted on the $y-$axis.  Within each chart, each bar represents 1000 trials, and each trial is $\epsilon$-noisy gradient descent run for 5000 steps.

From left to right, the charts are of step size $\tau$ equal to $0.001, 0.01, 0.02, 0.04,$ and $0.06$.  

We find that for all but the last value of $\tau$, there is for each choice of step size a distinct range of values of $\epsilon$ for which the ratio $r$ drops, and in some cases dramatically so.  For $\tau = 0.01, 0.02$, and $0.04$, there are ranges of $\epsilon$ for which $r$ drops to nearly 0.  Interestingly, this range is largest for $\tau = 0.02$, it appears to increase as $\tau$ increases until about $0.02$, and then decrease after that.  In the figure this appears as a gap that opens and gets wider till the middle chart, and then gets smaller.  It is also interesting that this gap shifts rightward as $\tau$ increases.  We leave the study of that phenomenon to future work.

Lastly, we note that for $\tau=0.06$, no gap appears.  This is not surprising, as in Section \ref{jitterthy} we noted that for $\tau$ larger than approximately $w/g$, discrete gradient descent in this landscape becomes more like jumping between wells and less like gradient flow.  So we did not expect that $\epsilon$ jitter would have a similar effect on discrete gradient descent for large $\tau$.  

\section{Gradient descent on $g(x)$: Theory}

Having seen that noise can bias discrete gradient descent toward deeper wells, we consider whether it can also bias discrete gradient descent toward wider wells.  To study this, we experiment with noisy discrete gradient descent on the function $g(x) =  \left(\sin(\pi x) + \frac{\sin(2 \pi x)}{2}\right)^2$.  

\subsection{Gradient flow}

As above, we consider the behavior of gradient flow on the function $g(x)$ on an interval $[a,b]$, with a given starting point $p_0 \in [a,b]$.  For simplicity, we assume that $a,b$ are maxima of $g$, and that the interval $[a,b]$ contains the same number of shallow wells as deep wells.  

As before, with measure 0, $p_0$ will be a critical point of $g$.  Assuming that $p_0$ is not a critical point of $g$, there are two possibilities.  Either $p_0$ is in the basin of attraction of a wide well $W$, or it is in the basin of attraction of a narrow well $N$.  While neither $W$ or $N$ is convex, nonetheless, under gradient flow if $p_0$ is in the basin of attraction of $W$ the flow line will end in the unique minimum of $W$, and similarly for $N$.  

Therefore, if we initialize $p_0$ randomly, we should expect that the ratio of the probability that a local nonglobal minimum is found to the probability that a global minimum is reached is the ratio of the width of the basin of attraction of $W$ to the width of the basin of attraction of $N$.  

The width of the basin of attraction of $N$ is $2/3$, and the width of the basin of attraction of $W$ is $4/3$, so their ratio is 0.5.  So we expect that under uniform random initialization, a minimum of a narrow well will be found approximately 0.5 as often as a minimum of a wide well.

\subsection{Discrete gradient descent with $\epsilon$-jitter}\label{with jitter}

As before, we will consider discrete gradient descent with $\epsilon$-jitter.  

For step sizes in the regime $0 < \tau < w/g$, we expect that the addition of noise of order $\epsilon$ to the gradient descent algorithm will affect the relative probabilities of finding wide and narrow minima.  We expect that for at least some range of choices for $\epsilon$, the addition of noise will bias discrete noisy gradient descent toward wide wells, because at each step, the probability of the noise causing the gradient descent path to leave a well is higher for narrower wells.  In the next section, we will experimentally study the dependence of that phenomenon on $\epsilon$, and compare the strength of the effect in this setting to the strength of the effect on deep vs. shallow wells.

\section{Discrete gradient descent on $g(x)$: Computer experiments}

We now turn to some computer experiments for gradient descent on $g(x)$.  In all the experiments of this section, we run variants of gradient descent on the function $g(x)$ over the interval $[-7/3,7/3]$, and initialize from the uniform distribution on that interval.

\subsection{Gradient flow}

First, we approximate gradient flow by discrete gradient descent with small step size.  We find, as expected, that the ratio of the probabilities of arriving at the minimum of a narrow well versus a wide well is equal to the ratio of the basins of attraction.  

\begin{figure}[h]
\includegraphics[width=3.5in]{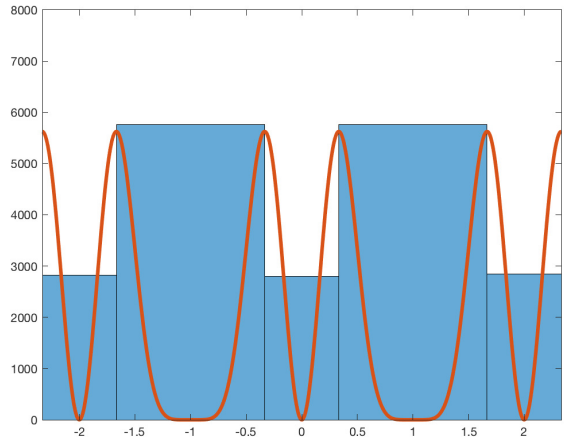}
\caption{In this experiment, we implement gradient descent with step size $\tau = 0.01$ on the function $g(x) =   \left(\sin(\pi x) + \frac{\sin(2 \pi x)}{2}\right)^2$ over the interval $[-7/3,7/3]$.  We randomly initialize from the uniform distribution on this interval, and run the experiment 20,000 times.  Above is a histogram which shows the number of times the process ended in each well.}
\end{figure}

Namely, the ratio we find in the experiment described above is 0.503, while the ratio between the widths of the basins of attraction of the two kinds of wells was computed above as 0.5.

\subsection{Discrete gradient descent with $\epsilon$-jitter}\label{with jitter}

In this section, we outline the results of some computer experiments carrying out discrete graident descent with $\epsilon$-jitter.  We find that indeed, for some values of $\epsilon$, adding noise to the gradient descent algorithm does bias the procedure toward finding the wider minima over the narrower ones.

In the first set of experiments, we implement gradient descent with step size $\tau = 0.01$ on the function $g(x) =   \left(\sin(\pi x) + \frac{\sin(2 \pi x)}{2}\right)^2$ over the interval $[-7/3,7/3]$.  We randomly initialize from the uniform distribution on this interval, and run the experiment 20,000 times.  We do this for several different values of noise $\epsilon$, and record the the number of times the process ended in each well.   

With $\epsilon = 0$, this experiment went similarly to the experiment approximating gradient flow, as expected.  With a relatively small amount of added noise however, the picture becomes dramatically different.

\begin{figure}[h]
\includegraphics[width=3.5in]{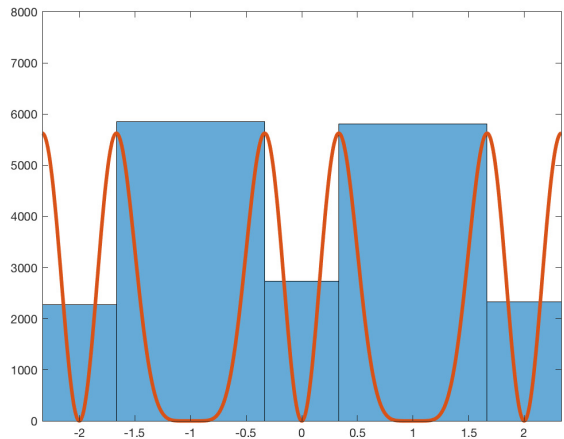}
\caption{This histogram shows the number of times the process ended in each well with $\epsilon = 0.15$.  The ratio of the probability of landing in a narrow well to the probability of landing in a wide well was .43.}
\end{figure}

With a small amount of noise added to the process at each step, the behavior of gradient descent on this function $g$ becomes noticeably different.  With $\epsilon = 0.15$, the probability of finding a shallow minimum decreases to approximately $0.43$.  

\begin{figure}[h]
\includegraphics[width=3.5in]{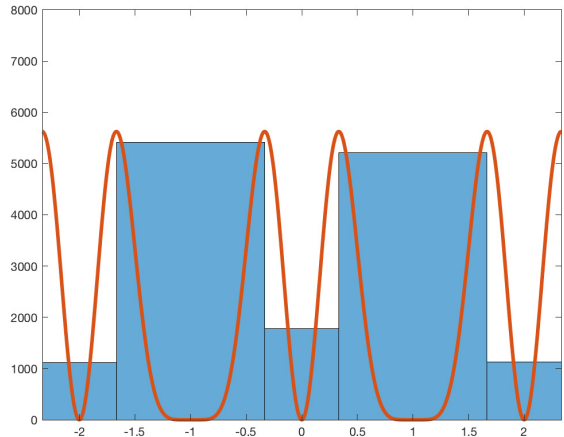}
\caption{This histogram shows the number of times the process ended in each well with $\epsilon = 0.25$.  The ratio of the probability of landing in a narrow well to the probability of landing in a wide well was 0.25.}
\end{figure}

The wide basin bias that $\epsilon$-jitter induces happens for a range of choices for $\epsilon$.  At $\epsilon = 0.25$, the effect is even more dramatic, and the probability of ending in a narrow well decreases to approximately $0.25$ the probability of ending in a wide well.

\begin{figure}[h]
\includegraphics[width=3.5in]{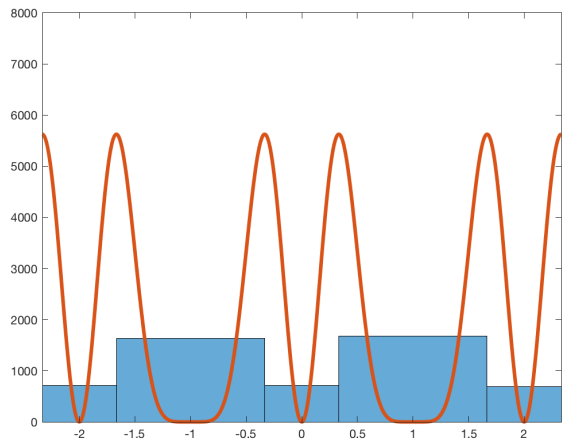}
\caption{This histogram shows the number of times the process ended in each well with $\epsilon = 0.5$.  The ratio of the probability of landing in a narrow well to the probability of landing in a wide well was 0.43.}
\end{figure}

Even for large values of $\epsilon$, when the noise is so large as to cause many trials to leave the interval of interest, the bias toward wide minima persists.  All the trials are initialized in the interval $[-7/3,7/3]$, but over the course of the 1,000 steps, most trials leave the interval.  Nonetheless, for those that remain in the interval, there is still a bias toward wide minima.

To address more systematically the question of the dependence of the behavior of discrete gradient descent with $\epsilon$-jitter on the step size $\tau$ and noise parameter $\epsilon$, we run experiments for a range of $\tau$ and $\epsilon$.  The results are collected in the following set of charts.  

\begin{figure}[h]\label{12}
\includegraphics[width=4in]{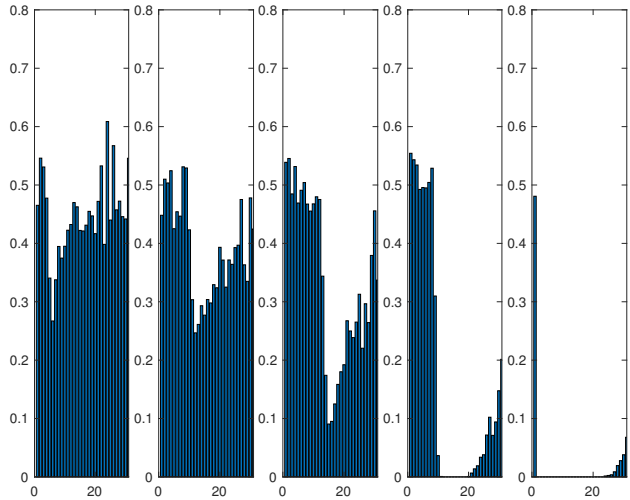}
\caption{These charts, from left to right, concern the behavior of noisy gradient descent on $g(x) =   \left(\sin(\pi x) + \frac{\sin(2 \pi x)}{2}\right)^2$, with step size $\tau = 0.001, 0.01, 0.02, 0.04$, and $0.06$.  Each chart shows, for that step size, how the ratio of the probability of ending in a deep vs a shallow well varies as the magnitude $\epsilon$ of the added noise varies from $0$ to $0.5$.}
\end{figure}

Similar to Figure \ref{6} from the previous section, each chart in Figure \ref{12} shows the ratio $r$ of the probability of finding a narrow minimum to the probability of finding a wide minimum when the step size $\tau$ is fixed, and the noise $\epsilon$ goes from 0 to $.5$.  Within each chart, each bar represents 1000 trials, and each trial is $\epsilon$-noisy gradient descent run for 5000 steps.  From left to right, the charts are of step size $\tau$ equal to $.001, .01, .02, .04,$ and $.06$.  

We find that for all but the last value of $\tau$, there is for each choice of step size a distinct range of values of $\epsilon$ for which the ratio $r$ drops, and in some cases dramatically so.  The ratio $r$ stays generally higher though than it did for the analogous experiment for the function $f(x)= \sin(\pi x) + \cos(2 \pi x)+2 $.  In this experiment, the only value of $\tau$ for which there is a range of $\epsilon$ for which $r$ drops nearly to 0 is $\tau = 0.04$.  And for $\tau = 0.04$ and $\epsilon = 0.2$, a value for which $r$ becomes very small, this is already in the regime where many trials end up outside the interval being studied.  

For smaller $\tau$, such as $\tau = 0.01, 0.02$, there are ranges of $\epsilon$ for which $r$ drops substantially, and as in the case studied in the previous section, we do again see a gap that appears around $\epsilon = 0.2$, and gets deeper and wider.  In this case, it does not get smaller again as $\tau$ increases.  It is interesting though that again, this gap shifts rightward as $\tau$ increases.  We leave the study of that phenomenon to future work.

In this experiment, in contrast to the one done in the case of deep vs. shallow wells, for $\tau=0.06$ the gap persists.  Since this step size is so large as to produce qualitatively different behavior from gradient flow even without noise, we did not make any prediction as to the effect of adding $\epsilon$-jitter in this setting, and understanding the role of $\epsilon$-jitter for large step size is a different problem than the one we are considering in this note.

\section{Discussion}

In this note, we ran computer experiments to explore the behavior of gradient descent and noisy discrete gradient descent.  We are interested in the role of noise in biasing the process toward and away from minima of various kinds.  To study this in a clean and simple setting, we focused on two landscapes.  Both are given by periodic functions.  In one, the graph of $f(x) = \sin(\pi x) + \cos(2 \pi x)+2$, there are two kinds of wells - deep and shallow wells.  In the other, the graph of $g(x)=  \left(\sin(\pi x) + \frac{\sin(2 \pi x)}{2}\right)^2$, there are again two kinds of wells.  Now the depths of the wells are the same but one is wide and the other narrow.

We find that in both cases, noise added to gradient descent can affect the probability of finding different minima.  In the first case, we saw that noise can bias the procedure toward finding deeper wells.  In the second case, we saw that noise can bias the procedure toward finding wider wells.  In both cases, the strength of the effect, and even whether it happens, is very sensitive to the step size $\tau$ and the magnitude of the noise $\epsilon$.  

We found that in both cases, when $\tau$ is small, as $\tau$ increases, the strength of the effect grows with $\tau$, and becomes more robust in terms of the range of choices for $\epsilon$ for which the effect is significant.  For larger $\tau$, the two cases were different in terms of how the strength of the effect changed when $\tau$ grew.  We did notice that in both cases, as $\tau$ increased, the values of $\epsilon$ for which the effect was strongest also increased.  We do not have an explanation for this, and leave the study of this for future work.

In these experiments, we did observe that adding noise to discrete gradient descent can bias the procedure toward deeper and wider wells, as is often suggested.  This mechanism may sometimes play a role in applications of gradient descent.

\end{document}